\documentclass{amsart}
\usepackage{amssymb}
\usepackage{braket}
\usepackage{mathrsfs}
\usepackage{ifthen}
\usepackage{here}
\usepackage{todonotes}
\usepackage{tikz}
\usetikzlibrary{patterns}
\usepackage{comment} 
\usepackage{mleftright}
\usepackage{hyperref}
\usepackage{cleveref}
 \usepackage[all]{xy}
 \usepackage{amscd}
 \usepackage{amsrefs}
 \usepackage{color}
 \usepackage{enumitem}
\newlist{steps}{enumerate}{1}
\setlist[steps, 1]{label = Step \arabic*:}
\usepackage[abs]{overpic}
\usepackage[latin1]{inputenc}
\usepackage{tikz-cd}

\makeatletter
\DeclareRobustCommand\widecheck[1]{{\mathpalette\@widecheck{#1}}}
\def\@widecheck#1#2{%
   \setbox\z@\hbox{\m@th$#1#2$}%
   \setbox\tw@\hbox{\m@th$#1%
      {%
         \vrule\@width\z@\@height\ht\z@
         \vrule\@height\z@\@width\wd\z@}$}%
   \dp\tw@-\ht\z@
   \@tempdima\ht\z@ \advance\@tempdima2\ht\tw@ \divide\@tempdima\thr@@
   \setbox\tw@\hbox{%
      \raise\@tempdima\hbox{\scalebox{1}[-1]{\lower\@tempdima\box\tw@}}}%
   {\ooalign{\box\tw@ \cr \box\z@}}}
\makeatother

\theoremstyle{plain}
\newtheorem{thm}{Theorem}[section]
\crefname{thm}{Theorem}{Theorems}
\Crefname{thm}{Theorem}{Theorems}

\crefname{prop}{Proposition}{Propositions}
\Crefname{prop}{Proposition}{Propositions}

\crefname{lem}{Lemma}{Lemmas}
\Crefname{lem}{Lemma}{Lemmas}
\newtheorem{cor}[thm]{Corollary}
\crefname{cor}{Corollary}{Corollaries}
\Crefname{cor}{Corollary}{Corollaries}

\crefname{claim}{Claim}{Claims}
\Crefname{claim}{Claim}{Claims}

\crefname{property}{Property}{Properties}
\Crefname{property}{Property}{Properties}

\crefname{problem}{Problem}{Problems}
\Crefname{problem}{Problem}{Problems}

\crefname{conjecture}{Conjecture}{Conjecture}
\Crefname{conjecture}{Conjecture}{Conjecture}

\theoremstyle{definition}
\newtheorem{defn}[thm]{Definition}
\crefname{defn}{Definition}{Definitions}
\Crefname{defn}{Definition}{Definitions}

\crefname{notation}{Notation}{Notations}
\Crefname{notation}{Notation}{Notations}

\crefname{convention}{Convention}{Conventions}
\Crefname{convention}{Convention}{Conventions}

\crefname{cond}{Condition}{Conditions}
\Crefname{cond}{Condition}{Conditions}

\crefname{assum}{Assumption}{Assumptions}
\Crefname{assum}{Assumption}{Assumptions}

\crefname{conj}{Conjecture}{Conjectures}
\Crefname{conj}{Conjecture}{Conjectures}

\crefname{claim1}{Claim}{Claims}
\Crefname{claim1}{Claim}{Claims}
\Crefname{ques}{Question}{Question}
\newtheorem{ques}[thm]{Question}
\crefname{que}{Question}{Question}
\Crefname{que}{Question}{Question}

\theoremstyle{remark}
\newtheorem{rem}[thm]{Remark}
\crefname{rem}{Remark}{Remarks}
\Crefname{rem}{Remark}{Remarks}
\newtheorem{ex}[thm]{Example}
\crefname{ex}{Example}{Examples}
\Crefname{ex}{Example}{Examples}

\crefname{section}{Section}{Sections}
\Crefname{section}{Section}{Sections}
\crefname{subsection}{Subsection}{Subsections}
\Crefname{subsection}{Subsection}{Subsections}
\crefname{figure}{Figure}{Figures}
\Crefname{figure}{Figure}{Figures}

\newtheorem*{acknowledgement}{Acknowledgement}

\newcommand{\Z}{\mathbb{Z}}

\newcommand{\Q}{\mathbb{Q}}

\newcommand{\Diff}{\mathrm{Diff}}
\newcommand{\Homeo}{\mathrm{Homeo}}

\newcommand{\del}{\partial}

\newcommand{\id}{\mathrm{id}}

\newcommand{\C}{\mathbb{C}}

\newcommand{\R}{\mathbb R}

\newcommand{\F}{\mathbb{F}_2}

\def\id{\operatorname{Id}}

\newcommand{\mbar}[1]{{\ooalign{\hfil#1\hfil\crcr\raise.167ex\hbox{--}}}}

     \RequirePackage{rotating}                   
    \def\HMt{%
       \setbox0=\hbox{$\widehat{\mathit{HM}}$}
       \setbox1=\hbox{$\mathit{HM}$}
       \dimen0=1.1\ht0
       \advance\dimen0 by 1.17\ht1
       \smash{\mskip2mu\raise\dimen0\rlap{%
          \begin{turn}{180}
              {$\widehat{\phantom{\mathit{HM}}}$}
           \end{turn}} \mskip-2mu    
                \mathit{HM}
                    }{\vphantom{\widehat{\mathit{HM}}}}{}}

\title[Exotic embeddings and stabilizations]{Exotic codimension-1 submanifolds in 4-manifolds and stabilizations}

\author{Hokuto Konno}
\address{Graduate School of Mathematical Sciences, the University of Tokyo, 3-8-1 Komaba, Meguro, Tokyo 153-8914, Japan \\and\\
RIKEN iTHEMS, Wako, Saitama 351-0198, Japan}
\email{konno@ms.u-tokyo.ac.jp}

\author{Anubhav Mukherjee}
\address{SL-Math(MSRI), 17 Gauss Way, Berkeley, CA 94720, USA}
\email{anubhavmaths@gmail.com}

\author{Masaki Taniguchi}
\address{RIKEN, 2-1 Hirosawa, Wako, Saitama 351-0198, Japan}
\email{masaki.taniguchi@riken.jp}

\date{}
\begin{document}

\begin{abstract}
In a small simply-connected closed 4-manifold, we construct infinitely many pairs of exotic codimension-$1$ submanifolds with diffeomorphic complements that remain exotic after any number of stabilizations by $ S^2\times S^2$.
We also give new constructions of exotic embeddings of 3-spheres in 4-manifolds with diffeomorphic complements.

\end{abstract}

\maketitle


\section{Introduction}

It is one of the key topics of the study in 4-manifold topology
to understand exotic phenomena, i.e. those properties that are true in the topological category but not in the smooth category. The main three exotic phenomena that have been getting lots of attention in recent times are existence of exotic smooth structures in 4-manifolds, existence of exotic surfaces in 4-manifolds and existence of exotic diffeomorphisms in 4-manifolds. 
Motivated from the 4-dimensional smooth Schoenflies conjecture, we recently started studying exotic codimension-1 submanifolds in 4-manifolds \cite{IKMT22}.
Here let us first clarify the term ``exotic codimension-1 submanifolds" of a 4-manifold and ``exotic embeddings of 3-manifolds" into a 4-manifold in this paper:

\begin{defn}
\label{def: exo submfd}
We say that two codimension-1 smooth submanifolds $Y_1$ and $Y_2$ in a smooth oriented 4-manifold $X$ (resp. two smooth embeddings $i_1,i_2 :Y\to X$) are {\it exotic} if \begin{itemize}
    \item[(i)] there is a topological ambient isotopy $H_t: X\times [0,1] \to X$ such that $H_0= \operatorname{Id}$ and $H_1(Y_1)=Y_2$ (resp. $H_1\circ i_1= i_2$),
    \item[(ii)] there is no such smooth isotopy,
    \item[(iii)] the complements of $Y_1$ and $Y_2$ are diffeomorphic, i.e there exists a diffeomorphism $f:X\to X$ such that $f(Y_1)= Y_2$ (resp. $f\circ i_1 = i_2$). 
\end{itemize} 
\end{defn}

Exotic codimension-1 submanifolds are closely related to interesting topics of diffeomorphisms of a 4-manifold.
For instance, the existence of exotic codimension-1 submanifolds leads to a counterexample to the $\pi_0$-version of the Smale conjecture in the following sense: 
If there is a pair of exotic codimension-1 submanifolds in $S^4$ whose diffeomorphism type is $Y$ so that the mapping class group of $Y$ is trivial, then it is easy to provide a non-trivial element of the kernel of
\[
\pi_0(\Diff (S^4)) \to \pi_0(\Homeo (S^4)). 
\]
This gives a way to disprove the 4-dimensional Smale conjecture different from the one by Watanabe~\cite{Wa18}.

\begin{rem}
\label{mcg equiv}
If the mapping class group of $Y$ is trivial, then exotic embeddings yield exotic submanifolds, and vice versa.
\end{rem}

\begin{rem}
The notion obtained by dropping the condition (iii) of \cref{def: exo submfd} might be called  ``weakly exotic submanifolds/embeddings".
It is relatively easy to construct such exotic codimension-1 submanifolds using (exotic) corks, as described in \cite[Introduction]{IKMT22}.
Corks yield also examples of this weaker version of exotic embeddings.
Indeed, let $(W,\tau)$ be a cork, where $\tau$ is a diffeomorphism of $\del W$.
Then it is easy to see that the double $W\cup_{\del W}(-W)$ admits weakly exotic embeddings of $\del W$ by considering two embeddings of $\del W$ that correspond to $\id_{\del W}$ and $\tau$.
However, this cork example does not give exotic embeddings in the sense of \cref{def: exo submfd}.
\end{rem}

We shall study the existence of exotic codimension-1 submanifolds, but also consider behavior of them under stabilizations, i.e. connected sum with $S^2\times S^2$'s, or more general 4-manifolds.
In the context of exotica  (in the category of orientable 4-manifolds), one remarkable discovery by Wall~\cite{wall} in 1960's says that 
exotic smooth structures disappear under sufficiently many stabilizations by $S^2\times S^2$.
For surfaces in a 4-manifold and diffeomorphisms of a 4-manifold,
analogous results were later established by Perron~\cite{P86}, and Quinn~\cite{Q86} (combined with a result of Kreck~\cite{Kreck79}).
One can similarly expect that such a principle should be true for exotic 3-manifolds in 4-manifolds.
Such a principle holds for previously known examples of exotic embeddings of 3-manifolds since those examples were constructed as images of exotic diffeomorphisms, and those diffeomorphisms are smoothly isotopic to the identity after stabilizations, as described.  

In order to state our first result, we use the Heegaard Floer tau-invariant \cite{OS03} 
$
\tau : \mathcal{C} \to \Z
$, where $\mathcal{C}$ denotes the knot concordance group.
Let $\mathrm{sign} : \Z\setminus\{0\} \to \{1, -1\}$ be the sign function.
The first main theorem of this paper is:

\begin{thm}
\label{cor submanifold} Let $K$ be a knot in $S^3$ with $\tau(K)\neq 0$ and let $n>0$. 
If the mapping class group of the 3-manifold $S^3_{\mathrm{sign}(\tau(K))/n}(K)$ is trivial, then 
\[
X_n := \begin{cases} \#_2 S^2 \times S^2 ,\text { $n$ is even } \\
\#_3 (\C P^2 \# (-\C P^2)) ,\text { $n$ is odd }
    \end{cases}
\] 
contains exotic codimension-1 submanifolds $Y_1$ and $Y_2$ whose diffeomorphism types are $S^3_{\mathrm{sign}(\tau(K))/n}(K)$, and which survives after a connected sum by any connected 4-manifold $M$, where $M$ is attached to a point away from $Y_1 \cup Y_2$. In particluar, we can choose $M= \#_m S^2\times S^2$ for any $m>0$.
\end{thm}

\begin{ex}
\label{ex: hyperbolic}
Many examples of knots satisfying the assumptions of \cref{cor submanifold} can be found as follows: 
Take a hyperbolic knot $K$ which has trivial symmetry group, i.e. the isometry group of the hyperbolic knot complement is trivial. For large $ n$, the mapping class group of $S^3_{1/n}(K)$ becomes trivial. (Here, we used the generalized Smale conjecture for hyperbolic 3-manifolds proven by Gabai \cite{Ga01} together with Thurston's hyperbolic Dehn filling theorem \cite{Thu}.) For example, using Snappy \cite{SnapPy}, we can ensure that one may take a hyperbolic knot $K=10_{149}$ with trivial symmetry group and with $\tau(K)=2$. 
\end{ex}

The authors do not know whether there exists a pair of codimension-1 submanifolds in smaller 4-manifolds than $\#_2 S^2\times S^2$ and $\#_3  (\C P^2 \# (-\C P^2))$. So we can ask the following question:

\begin{ques}
Does there exist a closed 4-manifold smaller than one in \cref{cor submanifold} that support exotic codimension-1 submanifolds?
More concretely, do $S^4$, $S^2\times S^2$, or $\C P^2\#(-\C P^2)$ contain exotic codimension-1 submanifolds? 
\end{ques}




\cref{cor submanifold} is derived from the following \lcnamecref{non-spin tau} on exotic embeddings, together with \cref{mcg equiv}:

\begin{thm}
\label{non-spin tau}
Let $K$ be a knot in $S^3$ with $\tau(K)\neq 0$ and let $n>0$. 
Then, there are two embeddings $i_1$ and $i_2$ from  $Y=S^3_{\mathrm{sign}(\tau(K))/n}(K)$ into 
\[
X_n := \begin{cases} \#_2 S^2 \times S^2 ,\text { $n$ is even } \\
\#_3 (\C P^2 \# (-\C P^2)) ,\text { $n$ is odd }
    \end{cases}
\] 
that are exotic. Moreover, they remain exotic in $X_n \# M$ for any closed smooth connected 4-manifold $M$, where $M$ is attached to a point away from $i_1(Y) \cup i_2(Y)$. In particular, they remain smoothly non-isotopic in $X_n \# (\#_m S^2\times S^2)$ for any $m>0$. 
\end{thm}

For example, we can take the positive torus knot $K=T(2,3)$. Then we have $S^3_{1/n} (K)\cong -\Sigma (2,3,6n-1)$. Moreover, we can take $K$ to be any alternating knot with non-zero knot signature or any non-slice quasi-positive knot. 
\begin{rem}
In recent work \cite{kang}, Kang showed the existence of a cork $W$ where the involution on the boundary $Y$ does not extend over the 4-manifold $W\# S^2\times S^2$ as a diffeomorphism. So if $X= W\cup_Y (-W)$, then we can find two embeddings of $Y$, i.e the identity map and the involution of $Y$ in $X$ that are topologically but not smoothly isotopic and remains exotic after connected sum a copy of $S^2\times S^2$ in the both components of the complement. In our language this is the first example of weakly exotic embeddings that survives stabilizations on both sides of the 3-manifolds. So we can ask the following question:
\end{rem}
\begin{ques}
Does there exist an example of exotic codimension-1 submanifolds/embeddings that survives after stabilizations on both sides of the complements? 
\end{ques}


 
 \begin{rem}
 Recently Lin \cite{Lin20} showed that the existence of an exotic diffeomorphism that survives a single stabilization and later Lin and the second author \cite{LM21} showed that there are exotic surfaces that survive after a single stabilization.
 \end{rem}
 
 As more concrete examples of exotic 3-manifolds, we can treat a class of homology 3-spheres bounded by definite 4-manifolds including $S^3$. 
 
\begin{thm}\label{contractible}
Let $Y$ be a homology 3-sphere with $d(Y)=0$ that bounds a compact simply-connected definite smooth 4-manifold, where $d$ denotes the Heegaard Floer d-invariant.  
Then there is a closed smooth 4-manifold $X$ containing a pair of exotic embeddings of $Y$.
\end{thm}

Since any Seifert homology 3-sphere bounds a simply-connected compact definite 4-manifold (for example, see \cite[Proposition 4.2]{HK18}), \cref{contractible} implies:

\begin{cor}\label{Seifert}
Every Seifert fibered homology 3-spheres with $d(Y)=0$ admits a pair of exotic embeddings into some closed smooth 4-manifold. 
\end{cor}


We also construct new examples of exotic 3-spheres in many 4-manifolds.

\begin{thm}
\label{thm: exotic emb spheres}
Let $X$ be a 4-manifold of the form 
\[
K3\# (S^2\times S^2)\# M \text{ or } 2\C P^2\#(-10\C P^2)\#M,
\]
where $M$ is any closed smooth 4-manifold. Then there exists a pair of exotic codimenion-1 submanifolds diffeomorphic to $S^3$.
\end{thm}
 
 As an immediate corollary of \cref{thm: exotic emb spheres} and Wall's theorem~\cite{wall}, we have the following result on exotic embeddings of 3-spheres. 
 
 \begin{cor}
Let $X$ be a simply-connected closed 4-manifold. Then there exists an integer $n>0$ such that, for any $k>n$,  $X\# k S^2\times S^2$ contains a pair of exotic codimenion-1 submanifolds diffeomorphic to $S^3$.
\end{cor}

\begin{rem}Gabai--Budney \cite{BG19} and 
Watanabe \cite{Wa20} have constructed examples of 3-balls in $S^1\times B^3$ and some general integer homology spheres in 4-manifolds that are not smoothly isotopic and their complements are diffeomorphic as they are related by some diffeomorphisms of ambient 4-manifolds. However it is not clear if their examples are topologically isotopic or not.
(Regarding the problem whether these examples are topologically standard or not, there is a paper \cite{Ma86} showing the topological light bulb theorem in dimension 4.)
In \cite{IKMT22}, Iida and the authors constructed the example of exotic 3-spheres in some 4-maniofolds, and to the authors best knowledge this is the first known example. In an upcoming paper entitled ``Exotic phenomena in dimension four: diffeomorphism groups and embedding spaces", Auckly and Ruberman also construct such examples. 
\cref{thm: exotic emb spheres} gives some new constructions of exotic 3-spheres in 4-manifold which are different from Auckly--Ruberman or our previous examples. 
\end{rem}

\begin{acknowledgement}
The authors wish to thank R. \.{I}nan{\c{c}} Baykur for pointing out an error of the proof of \cref{cor submanifold} in an earlier version of this paper. 
We also thank Ian Agol to tell us a fact on mapping class groups used in \cref{ex: hyperbolic} in \cite{ianagol}.
The authors would like to thank to Dave Auckly, John Etnyre, Michael Freedman, Kyle Hayden, Nobuo Iida, Tye Lidman, Jianfeng Lin, Ciprian Manolescu, Danny Ruberman and Tadayuki Watanabe for various discussions and giving helpful comments on our first draft.

This material is based upon work supported by the National Science Foundation under Grant No. DMS-1928930, while the authors were in residence at the Simons Laufer Mathematical Science Institute (previously known as MSRI) Berkeley, California, during fall 2022 semester. 

In addition, Hokuto Konno was partially supported by JSPS KAKENHI Grant Numbers 17H06461, 19K23412, and 21K13785,
Anubhav Mukherjee was partially supported by NSF grant DMS-1906414 and Masaki Taniguchi was partially supported by JSPS KAKENHI Grant Number 20K22319, 22K13921 and RIKEN iTHEMS Program.
\end{acknowledgement}

\section{Proof of main results}

\subsection{Non-smoothable homeomorphisms}

Let us write the Fr\o yshov invariant by $\delta (Y) \in \Z$
which is defined using Seiberg--Witten Floer homotopy type \cite{Man03} and which has several different formulations \cite{Fr96, Fr10, KM07}. For example, see \cite[Subsection 2.2]{KT20}.  
The main ingredient in this paper is as follows: 

\begin{thm}\label{main theo}
Let $X$ be a compact simply-connected oriented smooth 4-manifold with $\sigma(X)\leq0$, $b^+(X)=1$ and with boundary $Y$ which is an integral homology 3-sphere.
Suppose that $X$ and $Y$ satisfy at least one of the following conditions:
\begin{enumerate}
\item $\delta(Y) \leq 0$ and $\sigma(X)<-8$.
\item $\delta(Y) < 0$, and in addition $\sigma(X)<0$ if $X$ is non-spin.
\end{enumerate}
Then there is no diffeomorphism $g :X \to X$ which reverses an orientation of $H^+(X)$ and fixes boundary pointwise. Moreover, there is a homeomorphism $h: X \to  X$ which reverses an orientation of $H^+(X)$ and fixes boundary pointwise. 

\end{thm}

We mainly use the techniques developed in \cite{KT20} to prove \cref{main theo}.  

\begin{proof}[Proof of \cref{main theo}]
The first fact, i.e. existence of no such diffeomorphism $g$ follows from \cite[Theorem 1.1]{KT20}. 
In \cite{KT20}, under the assumption of \cref{main theo}, we constructed a relative homeomorphism $h$ (i.e. homeomorphism that is the identity on the boundary) of $X$ which is not smoothable as a relative diffeomorphism i.e. not topologically isotopic to any relative diffeomorphism.

Since $h$ is the main ingredient of this paper, we give a sketch of the construction in the case of \cref{main theo} (1) for the reader's convenience.
Suppose that $\sigma(X)<-8$. 
Let $-E_{8}$ denote the negative-definite topological $E_{8}$-manifold and  $-\C P^{2}_{\mathrm{fake}}$ denote the fake $-\C P^{2}$.
Then $X$ is homeomorphic to
\begin{align}
\label{eq: nonspin connected sum ho}
S^{2} \times S^{2} \# n(-\C P^{2}) \# (-E_{8}) \#(-\C P^{2}_{\mathrm{fake}})\#W
\end{align}
or
\begin{align}
\label{eq: nonspin connected sum}
S^{2} \times S^{2} \# n(-\C P^{2}) \# (-E_{8}) \#W, 
\end{align}
where $n \geq 0$ and $W$ is a compact oriented contractible topological $4$-manifold bounded by $Y$.
Let $f_{0}$ be an orientation-preserving self-diffeomorphism of $S^2\times S^2$ which reverses orientation of $H^+(S^2\times S^2)$ and has a fixed 4-disk, where $H^+(X)$ is a maximal positive definite subspace of $H^2(X; \R)$ with respect to the intersection form of a given 4-manifold $X$. 
Such $f_{0}$ can be easily made by deforming the componentwise complex conjugation of $\C P^{1} \times \C P^{1} = S^{2} \times S^{2}$ to have a pointwise fixed 4-disk.

We extend $f_0$ as a homeomorphism of $X$ by the identity map on the other connected sum factors in \eqref{eq: nonspin connected sum ho} or \eqref{eq: nonspin connected sum}. We call the extension  $h$.
Note that the restriction of $h$ onto the boundary is the identity map. 
Thus $h$ satisfies desired condition. 
\end{proof}

\subsection{Proof of the main results}

Before proving \cref{non-spin tau}, let us recall the following two facts: 
\begin{rem}
The following fact is pointed out in \cite[Remark1.1]{LRS18}. 
The isomorphism between the monopole Floer homology and the Heegaard Floer homology \cite{KLTI}, \cite{KLTII}, \cite{KLTIII}, \cite{KLTIV}, \cite{KLTV}, alternatively, the work of Colin, Ghiggini, and Honda \cite{CGHI} \cite{CGHII} \cite{CGHIII}  and Taubes \cite{Ta10}, 
 (for the comparison of $\Q$-gradings: see \cite{RG18}, \cite{CD13} and \cite{HR17})
implies 
\[
\frac{1}{2} d ( Y, \mathfrak{t}, \F  ) =   -h ( Y, \mathfrak{t},  \F   ), 
\]
where $d ( Y, \mathfrak{t}, \F  )$ is the correction term of Heegaard Floer homology defined over $\F$-coefficients and $h ( Y, \mathfrak{t},  \F   )$ is the monopole Fr\o yshov invariant defined over $\F$-coefficients. 
On the other hand, in \cite{LM18}, it is proven that 
\[
-h ( Y, \mathfrak{t},  \F   ) = \delta (Y, \mathfrak{t}). 
\]
Summarizing the results above, we have 
\begin{align}\label{d=delta}
\frac{1}{2} d ( Y, \mathfrak{t}, \F  ) =  \delta (Y, \mathfrak{t}).
\end{align}
In this paper, we use the Heegaard Floer $d$-invariant over $\Z_2$-cofficients which we will denote by $d(Y)$ for an oriented homology 3-sphere $Y$ with unique spin$^c$ structure.
\end{rem}

\begin{rem}
We remark that there is the following relation between the correction term $d$ and Heegaard Floer tau-invariant $\tau$: 
\begin{align}\label{surgery formula of d}
\text{ if }\tau(K) >0 \text{ then } d(S_{1/n}^3(K))<0
\end{align}
for any positive integer $n$. 
First,  \eqref{surgery formula of d} is proven for $n=1$ essentially in \cite{HW16}. (See also \cite{Sa18}.) Then \eqref{surgery formula of d} for general $n$ follows from the fact that
\begin{align*}
    d(S_1^3(K)) = d(S_{1/n}^3(K)),
\end{align*}
which is proven in \cite[Proposition 1.6]{NW15}. 
\end{rem}

Now we prove the first main theorem of this paper, examples of exotic codimenion-1 embeddings/submanifolds in small 4-manifolds that survive after any stabilizations:

\begin{proof}[Proof of \cref{non-spin tau}]
First we prove for the case when $\tau(K)>0$.
From the remarks \eqref{d=delta} and \eqref{surgery formula of d}, we have 
\begin{align}
\label{eq: delta in the proof of main}
\delta (S_{1/n} ^3(K))= \delta (S_1^3(K))<0. 
\end{align}
Let $W_{1/n}(K)$ be the trace of the surgery with respect to the left Kirby diagram in Figure~1 with boundary $Y=S^3_{1/n}(K)$.
The intersection form of $W_{1/n}(K)$ is isomorphic to  $\mathrm{diag}(1,-1)$ for $n$ odd, and isomorhic  to $\begin{pmatrix}
0 & 1 \\
1 & 0
\end{pmatrix}$ for $n$ even.
Set 
\[
X_n':= W_{1/n}(K)\cup_{S^3_{1/n}(K)} (-W_{1/n}(K)).
\]
Then, the handle slidings along the red 0-framed 2-handles shown in Figure 1 proves that there is a diffeomorphism
\[
X_n' \cong  \begin{cases} \#_2 S^2 \times S^2 \quad \text {if $n$ is even, } \\
\#_2 (\C P^2 \# (-\C P^2))  \quad \text { if $n$ is odd. }
\end{cases}
\]

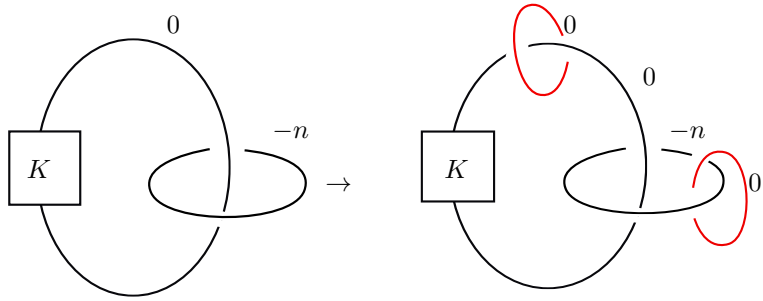
\begin{figure}[htbp]\label{Kirby}
	\begin{center}
  \begin{overpic}[scale=0.6,  tics=20]{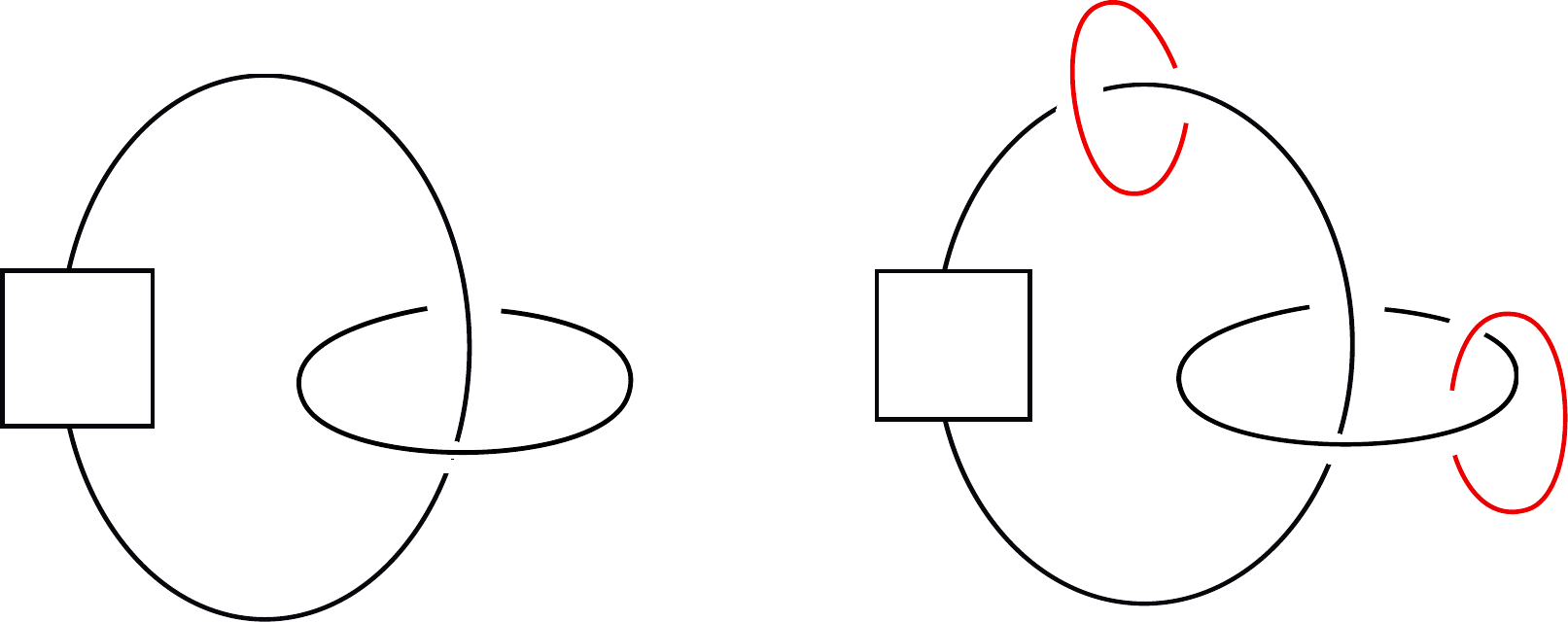}
  \put(7,45){$K$}
  \put(165,45){$K$}
  \put(60,100){$0$}
  \put(100,60){$-n$}
  \put(120,40){$\rightarrow$}
  \put(240,80){$0$}
  \put(250,60){$-n$}
  \put(280,40){$0$}
  \put(210,100){$0$}

\end{overpic}

\caption{Taking a double of $W_{1/n}(K)$}
\label{cob}
\end{center}
\end{figure}

Now we set
\[
X_n= \begin{cases} X_n'= W_{1/n}(K)\cup_{Y} (-W_{1/n}(K)) \quad \text {if $n$ is even, } \\
(-\C P^2)\# W_{1/n}(K)\cup_{Y} (-W_{1/n}(K)) \# \C P^2  \quad \text { if $n$ is odd. }
\end{cases}  
\]
Let $i_1 : Y \to X_n$ be a smooth embedding obtained as the inclusion.
Now, for an even $n$ (resp. odd $n$), 
let $h$ be the homeomorphism of $W_{1/n}(K)$ (resp. $(-\C P^2) \# W_{1/n}(K)$) that fixes the boundary pointwise, given in \cref{main theo}.
We define a homeomorphism $h'$ of $X_n$ by 
\[
h':= h \cup_{\id_Y} (-h) : X_n \to X_n.
\]


It is easy to see by an elementary argument (or by Wall's theorem \cite{Wa64}) that
any automorphism of the intersection form of $X_n$ can be realized by an orientation-preserving diffeomorphism.
In particular, there is an orientation-preserving diffeomorphism $f$ of $X_n'$ such that 
\[
f_* = h'_* : H_2(X_n;\Z) \to H_2(X_n;\Z).
\]
Define a smooth embedding $i_2: Y \to X_n$ by
\[
i_2=f \circ i_1 : Y \to X_n. 
\]
A result of Perron~\cite{P86} and Quinn~\cite{Q86} implies that there is a topological isotopy from $h'$ to $f$.
Also, $h'$ fixes $Y$ pointwise by construction.
It follows from these that $i_1$ and $i_2$ are mutually topologically isotopic.

However, $i_1$ and $i_2$ are not smoothly isotopic to each other. Indeed, if so, for an even $n$ (resp. odd $n$), by composing an ambient isotopy between $i_1$ and $i_2$ with $f$, we obtain a diffeomorphism $g$ of $W_{1/n}(K)$ (resp. of $(-\C P^2)\#W_{1/n}(K)$) that preserves $Y$ pointwise and that reverses orientation of $H^+(W_{1/n}(K))$ (resp. $H^+((-\C P^2)\#W_{1/n}(K))$).
This combined with \eqref{eq: delta in the proof of main} contradicts \cref{main theo}. 

In the case of $\tau(K)<0$, we have $\tau(K^*)>0$ for the mirror image $K^\ast$ of $K$. Again, from the remarks \eqref{d=delta} and \eqref{surgery formula of d}, we have 
\[
\delta (S^3_{1/n}(K^*) )=\delta (S_1^3(K^*))<0. 
\]
Since $S^3_{1/n}(K^*) \cong -S^3_{-1/n}(K)$ and $\delta$ is a homomorphism, we have 
\[
\delta (S^3_{-1/n}(K) )=\delta (S_{-1}^3(K))>0. 
\]
Then, we can do the same discussion for the 4-manifold corresponding the 4-manifold $-W_{1/n}$ with boundary  in Figure 1.

Next, we take a connected sum with a connected 4-manifold $M$ attached to a point $x \in X_n$ away from $i_1(Y) \cup i_2(Y)$. Suppose that $i_1$ and $i_2$ are smoothly isotopic. Then, for an even $n$ (resp. odd $n$), after applying an ambient isotopy between $i_1$ and $i_2$, the point $x$ lies in one of $ W_{1/n}(K)$ or $-W_{1/n}(K)$ (resp. $(-\C P^2)\# W_{1/n}(K)$ or $(-W_{1/n}(K)) \#\C P^2$). 
For an even $n$ (resp. odd $n$), if $x \in  W_{1/n}(K)$ (resp. $x \in(-\C P^2)\# W_{1/n}(K))$, we can do the same argument as before using the opposite orientation of $-W_{1/n}(K)$ (resp. $(-W_{1/n}(K)) \#\C P^2$). If $x \in -W_{1/n}(K)$ (resp. $x \in (-W_{1/n}(K)) \#\C P^2$), we do the same argument as before using $ W_{1/n}(K)$ (resp. $(-\C P^2)\# W_{1/n}(K)$). This completes the proof. 
\end{proof}

We give the proof of \cref{contractible} which gives more examples of exotic codimension-1 embeddings. 

\begin{proof}[Proof of \cref{contractible}]
We know that $\delta(Y)=0$ from \eqref{d=delta}. 
Let $W$ be a simply-connected definite smooth compact 4-manifold bounded by a homology 3-sphere $Y$.
Without loss of generality, we can assume that $W$ is negative definite. Take a simply-connected smooth 4-manifold $W'$ bounded by $-Y$ and consider the 4-manifold 
\[
X: = \C P^2 \# (-10\C P^2)  \# W \cup_Y W'. 
\]
If necessary, we may assume that $W'$ contains $S^2 \times S^2$ as a summand of a smooth connected sum decomposition.
Now, we apply \cref{main theo}~(1) to $\C P^2 \# (-10\C P^2)  \# W$ to obtain a homeomorphism $h'$ that fixes the boundary pointwise and reverses orientation of $H^+$.
Let $h: X \to X$ be a homeomorphism defined as the extension of $h'$ by the identity of $W'$. 
As in the proof of \cref{non-spin tau},
combinining Wall's theorem with Quinn and Perron's theorem, we may take a diffeomorphism $g$ on $X$ such that $g$ is topologically isotopic to $h$.  Then we see the natural inclusion $i:Y\to X$ and $g\circ i:Y\to X$ are exotic by a similar argument in \cref{non-spin tau}. 
\end{proof}



Lastly, we prove the result on exotic embeddings of 3-spheres: 

\begin{proof}[Proof of \cref{thm: exotic emb spheres}]
We start with the case that $X$ is of the form $X=K3\# S^2\times S^2\# M$. Then, from the topological decomposition $K3 \cong \#_3 S^2\times S^2 \# (\#_2 (-E_8)) $, one can construct a homeomorphism $h'$ on $K3 \setminus D^4$ rel boundary so that $h'$ acts on $H^+ (K3 \setminus D^4)$ as an orientation reversing map. Then we define a homeomorphism $h$ on $K3\#_{S^3} S^2\times S^2\# M$ by extending $h'$ as the identity on $(S^2\times S^2 \setminus D^4) \# M$. 

Let $Y=S^3$ be a separating 3-sphere in the neck of the connected sum of $K3$ and $S^2\times S^2 \# M$. Now, using Wall's theorem and Quinn and Perron's result, we may take an orientation-preserving diffeomorphism $f$ of $X$ which is topologically isotopic to $h$. 

We claim that $Y$ and $f(Y)$ are exotic submanifolds.
First, $Y$ and $f(Y)$ are topologically isotopic since $f$ is topologically isotopic to $h$ and $h$ fixes $Y$ pointwise.
Suppose that $Y$ and $f(Y)$ are smoothly isotopic.
Then, we obtain a diffeomorphism $g'$ on $K3$ that reverses orientation of $H^+(K3)$ noting that the mapping class group of $S^3$ is trivial.
This contradicts a result by Donaldson~\cite{Do90} (or \cref{main theo}~(1)).

Next, we consider the case that $X=2\C P^2\#(-10\C P^2)\#M$. 
We apply \cref{main theo} (1) to $
\C P^2\# (-10\C P^2 \setminus D^4)$  and obtain a homeomorphism $h'$ on $
\C P^2\# (-10\C P^2 \setminus D^4)$ rel boundary which cannot be isotopic to any diffeomorphism. We extend $h'$ to a homeomorphism $h$ on $X$ by putting the identity on the rest of $X$. We again use Wall's theorem and Quinn and Perron's theorem to obtain a diffeomorphism $g$ which is topologically isotopic to $h$.
Now we consider the 3-sphere $Y=S^3$ that separates $\C P^2\#(-10\C P^2)$ and $\C P ^2 \#M$ in $X$.

We claim that $Y$ and $ g(Y)$ are exotic submanifolds.
Again these are topologically isotopic by construction.
If they are smoothly isotopic,
we obtain a diffeomorphism of $\C P^2\# (-10\C P^2)$ that reverses orientation of $H^+$.
This contradicts a result by  Friedman--Morgan~\cite{FM88} (or \cref{main theo}~(1)).
This completes the proof. 
\end{proof}

\bibliographystyle{plain}
\bibliography{tex}

\begin{bibdiv}
\begin{biblist}

\bib{BG19}{article}{
      author={Budney, Ryan},
      author={Gabai, David},
       title={Knotted 3-balls in {$S^4$}},
        date={2019},
      eprint={arXiv:1912.09029},
}

\bib{HK18}{article}{
      author={Choe, Dong~Heon},
      author={Park, Kyungbae},
       title={On intersection forms of definite 4-manifolds bounded by a
  rational homology 3-sphere},
        date={2018},
        ISSN={0166-8641},
     journal={Topology Appl.},
      volume={238},
       pages={59\ndash 75},
         url={https://doi.org/10.1016/j.topol.2018.01.013},
      review={\MR{3775118}},
}

\bib{CGHIII}{article}{
      author={Colin, Vincent},
      author={Ghiggini, Paolo},
      author={Honda, Ko},
       title={The equivalence of heegaard floer homology and embedded contact
  homology {III}: from hat to plus},
        date={2012},
      eprint={arXiv:1208.1526},
}

\bib{CGHI}{article}{
      author={Colin, Vincent},
      author={Ghiggini, Paolo},
      author={Honda, Ko},
       title={The equivalence of {H}eegaard {F}loer homology and embedded
  contact homology via open book decompositions {I}},
        date={2012},
      eprint={arXiv:1208.1074},
}

\bib{CGHII}{article}{
      author={Colin, Vincent},
      author={Ghiggini, Paolo},
      author={Honda, Ko},
       title={The equivalence of {H}eegaard {F}loer homology and embedded
  contact homology via open book decompositions {II}},
        date={2012},
      eprint={arXiv:1208.1077},
}

\bib{CD13}{article}{
      author={Cristofaro-Gardiner, Daniel},
       title={The absolute gradings on embedded contact homology and
  {S}eiberg-{W}itten {F}loer cohomology},
        date={2013},
        ISSN={1472-2747},
     journal={Algebr. Geom. Topol.},
      volume={13},
      number={4},
       pages={2239\ndash 2260},
         url={https://doi.org/10.2140/agt.2013.13.2239},
      review={\MR{3073915}},
}

\bib{SnapPy}{misc}{
      author={Culler, Marc},
      author={Dunfield, Nathan~M.},
      author={Goerner, Matthias},
      author={Weeks, Jeffrey~R.},
       title={Snap{P}y, a computer program for studying the geometry and
  topology of $3$-manifolds},
         how={Available at \url{http://snappy.computop.org} (DD/MM/YYYY)},
}

\bib{Do90}{article}{
      author={Donaldson, S.~K.},
       title={Polynomial invariants for smooth four-manifolds},
        date={1990},
        ISSN={0040-9383},
     journal={Topology},
      volume={29},
      number={3},
       pages={257\ndash 315},
         url={https://doi.org/10.1016/0040-9383(90)90001-Z},
      review={\MR{1066174}},
}

\bib{FM88}{article}{
      author={Friedman, Robert},
      author={Morgan, John~W.},
       title={On the diffeomorphism types of certain algebraic surfaces. {I}},
        date={1988},
        ISSN={0022-040X},
     journal={J. Differential Geom.},
      volume={27},
      number={2},
       pages={297\ndash 369},
         url={http://projecteuclid.org/euclid.jdg/1214441784},
      review={\MR{925124}},
}

\bib{Fr96}{article}{
      author={Fr\o~yshov, Kim~A.},
       title={The {S}eiberg-{W}itten equations and four-manifolds with
  boundary},
        date={1996},
        ISSN={1073-2780},
     journal={Math. Res. Lett.},
      volume={3},
      number={3},
       pages={373\ndash 390},
         url={https://doi.org/10.4310/MRL.1996.v3.n3.a7},
      review={\MR{1397685}},
}

\bib{Fr10}{article}{
      author={Fr{\o}yshov, Kim~A.},
       title={Monopole {F}loer homology for rational homology 3-spheres},
        date={2010},
        ISSN={0012-7094},
     journal={Duke Math. J.},
      volume={155},
      number={3},
       pages={519\ndash 576},
         url={https://doi.org/10.1215/00127094-2010-060},
      review={\MR{2738582}},
}

\bib{Ga01}{article}{
      author={Gabai, David},
       title={The {S}male conjecture for hyperbolic 3-manifolds: {${\rm
  Isom}(M^3)\simeq{\rm Diff}(M^3)$}},
        date={2001},
        ISSN={0022-040X},
     journal={J. Differential Geom.},
      volume={58},
      number={1},
       pages={113\ndash 149},
         url={http://projecteuclid.org/euclid.jdg/1090348284},
      review={\MR{1895350}},
}

\bib{HW16}{article}{
      author={Hom, Jennifer},
      author={Wu, Zhongtao},
       title={Four-ball genus bounds and a refinement of the
  {O}zv\'{a}th-{S}zab\'{o} tau invariant},
        date={2016},
        ISSN={1527-5256},
     journal={J. Symplectic Geom.},
      volume={14},
      number={1},
       pages={305\ndash 323},
         url={https://doi.org/10.4310/JSG.2016.v14.n1.a12},
      review={\MR{3523259}},
}

\bib{ianagol}{article}{
      author={(https://mathoverflow.net/users/1345/ian agol), Ian~Agol},
       title={On trivial mapping class group of 3-manifolds},
         how={MathOverflow},
      eprint={https://mathoverflow.net/q/432729},
         url={https://mathoverflow.net/q/432729},
        note={URL:https://mathoverflow.net/q/432729 (version: 2022-10-21)},
}

\bib{HR17}{article}{
      author={Huang, Yang},
      author={Ramos, Vinicius G.~B.},
       title={An absolute grading on {H}eegaard {F}loer homology by homotopy
  classes of oriented 2-plane fields},
        date={2017},
        ISSN={1527-5256},
     journal={J. Symplectic Geom.},
      volume={15},
      number={1},
       pages={51\ndash 90},
         url={https://doi.org/10.4310/JSG.2017.v15.n1.a2},
      review={\MR{3652073}},
}

\bib{IKMT22}{article}{
      author={Iida, Nobuo},
      author={Konno, Hokuto},
      author={Mukherjee, Anubhav},
      author={Taniguchi, Masaki},
       title={Diffeomorphisms of 4-manifolds with boundary and exotic
  embeddings},
        date={2022},
      eprint={arXiv:2203.14878},
}

\bib{kang}{article}{
      author={Kang, Sungkyung},
       title={One stabilization is not enough for contractible 4-manifolds},
        date={2022},
         url={https://arxiv.org/abs/2210.07510},
}

\bib{KT20}{article}{
      author={Konno, Hokuto},
      author={Taniguchi, Masaki},
       title={The groups of diffeomorphisms and homeomorphisms of 4-manifolds
  with boundary},
        date={2022},
        ISSN={0001-8708},
     journal={Adv. Math.},
      volume={409},
       pages={Paper No. 108627},
         url={https://doi.org/10.1016/j.aim.2022.108627},
      review={\MR{4469073}},
}

\bib{Kreck79}{inproceedings}{
      author={Kreck, M.},
       title={Isotopy classes of diffeomorphisms of {$(k-1)$}-connected
  almost-parallelizable {$2k$}-manifolds},
        date={1979},
   booktitle={Algebraic topology, {A}arhus 1978 ({P}roc. {S}ympos., {U}niv.
  {A}arhus, {A}arhus, 1978)},
      series={Lecture Notes in Math.},
      volume={763},
   publisher={Springer, Berlin},
       pages={643\ndash 663},
      review={\MR{561244}},
}

\bib{KM07}{book}{
      author={Kronheimer, Peter},
      author={Mrowka, Tomasz},
       title={Monopoles and three-manifolds},
      series={New Mathematical Monographs},
   publisher={Cambridge University Press, Cambridge},
        date={2007},
      volume={10},
        ISBN={978-0-521-88022-0},
         url={https://doi.org/10.1017/CBO9780511543111},
      review={\MR{2388043}},
}

\bib{KLTIV}{article}{
      author={Kutluhan, Ca\u{g}atay},
      author={Lee, Yi-Jen},
      author={Taubes, Clifford},
       title={H{F}={HM}, {IV}: {T}he {S}ieberg-{W}itten {F}loer homology and
  ech correspondence},
        date={2020},
        ISSN={1465-3060},
     journal={Geom. Topol.},
      volume={24},
      number={7},
       pages={3219\ndash 3469},
         url={https://doi.org/10.2140/gt.2020.24.3219},
      review={\MR{4194308}},
}

\bib{KLTV}{article}{
      author={Kutluhan, Ca\u{g}atay},
      author={Lee, Yi-Jen},
      author={Taubes, Clifford~Henry},
       title={H{F}={HM}, {V}: {S}eiberg-{W}itten {F}loer homology and handle
  additions},
        date={2020},
        ISSN={1465-3060},
     journal={Geom. Topol.},
      volume={24},
      number={7},
       pages={3471\ndash 3748},
         url={https://doi.org/10.2140/gt.2020.24.3471},
      review={\MR{4194309}},
}

\bib{KLTI}{article}{
      author={Kutluhan, Ca\u{g}atay},
      author={Lee, Yi-Jen},
      author={Taubes, Clifford~Henry},
       title={{$\rm HF{=}HM$}, {I}: {H}eegaard {F}loer homology and
  {S}eiberg-{W}itten {F}loer homology},
        date={2020},
        ISSN={1465-3060},
     journal={Geom. Topol.},
      volume={24},
      number={6},
       pages={2829\ndash 2854},
         url={https://doi.org/10.2140/gt.2020.24.2829},
      review={\MR{4194305}},
}

\bib{KLTII}{article}{
      author={Kutluhan, Ca\u{g}atay},
      author={Lee, Yi-Jen},
      author={Taubes, Clifford~Henry},
       title={{$\rm HF{=}HM$}, {II}: {R}eeb orbits and holomorphic curves for
  the ech/{H}eegaard {F}loer correspondence},
        date={2020},
        ISSN={1465-3060},
     journal={Geom. Topol.},
      volume={24},
      number={6},
       pages={2855\ndash 3012},
         url={https://doi.org/10.2140/gt.2020.24.2855},
      review={\MR{4194306}},
}

\bib{KLTIII}{article}{
      author={Kutluhan, Ca\u{g}atay},
      author={Lee, Yi-Jen},
      author={Taubes, Clifford~Henry},
       title={{$\rm HF{=}HM$}, {III}: holomorphic curves and the differential
  for the ech/{H}eegaard {F}loer correspondence},
        date={2020},
        ISSN={1465-3060},
     journal={Geom. Topol.},
      volume={24},
      number={6},
       pages={3013\ndash 3218},
         url={https://doi.org/10.2140/gt.2020.24.3013},
      review={\MR{4194307}},
}

\bib{LM18}{article}{
      author={Lidman, Tye},
      author={Manolescu, Ciprian},
       title={The equivalence of two {S}eiberg-{W}itten {F}loer homologies},
        date={2018},
        ISSN={0303-1179},
     journal={Ast\'{e}risque},
      number={399},
       pages={vii+220},
      review={\MR{3818611}},
}

\bib{Lin20}{misc}{
      author={Lin, Jianfeng},
       title={Isotopy of the {D}ehn twist on ${K}3\#{K}3$ after a single
  stabilization},
        date={2020},
}

\bib{LM21}{article}{
      author={Lin, Jianfeng},
      author={Mukherjee, Anubhav},
       title={Family {B}auer--{F}uruta invariant, exotic surfaces and smale
  conjecture},
        date={2021},
      eprint={arXiv:2110.09686},
}

\bib{LRS18}{article}{
      author={Lin, Jianfeng},
      author={Ruberman, Daniel},
      author={Saveliev, Nikolai},
       title={{O}n the {F}r{\o}yshov invariant and monopole {L}efschetz
  number},
        date={2018},
      eprint={arXiv:1802.07704},
        note={to appear in Journal of Differential Geometry},
}

\bib{Man03}{article}{
      author={Manolescu, Ciprian},
       title={Seiberg-{W}itten-{F}loer stable homotopy type of three-manifolds
  with {$b_1=0$}},
        date={2003},
        ISSN={1465-3060},
     journal={Geom. Topol.},
      volume={7},
       pages={889\ndash 932},
         url={https://doi.org/10.2140/gt.2003.7.889},
      review={\MR{2026550}},
}

\bib{Ma86}{article}{
      author={Marumoto, Yoshihiko},
       title={On higher-dimensional light bulb theorem},
        date={1986},
        ISSN={0289-9051},
     journal={Kobe J. Math.},
      volume={3},
      number={1},
       pages={71\ndash 75},
      review={\MR{867805}},
}

\bib{NW15}{article}{
      author={Ni, Yi},
      author={Wu, Zhongtao},
       title={Cosmetic surgeries on knots in {$S^3$}},
        date={2015},
        ISSN={0075-4102},
     journal={J. Reine Angew. Math.},
      volume={706},
       pages={1\ndash 17},
         url={https://doi.org/10.1515/crelle-2013-0067},
      review={\MR{3393360}},
}

\bib{OS03}{article}{
      author={Ozsv\'{a}th, Peter},
      author={Szab\'{o}, Zolt\'{a}n},
       title={On the {F}loer homology of plumbed three-manifolds},
        date={2003},
        ISSN={1465-3060},
     journal={Geom. Topol.},
      volume={7},
       pages={185\ndash 224},
         url={https://doi.org/10.2140/gt.2003.7.185},
      review={\MR{1988284}},
}

\bib{P86}{article}{
      author={Perron, B.},
       title={Pseudo-isotopies et isotopies en dimension quatre dans la
  cat\'{e}gorie topologique},
        date={1986},
        ISSN={0040-9383},
     journal={Topology},
      volume={25},
      number={4},
       pages={381\ndash 397},
         url={https://doi.org/10.1016/0040-9383(86)90018-2},
      review={\MR{862426}},
}

\bib{Q86}{article}{
      author={Quinn, Frank},
       title={Isotopy of {$4$}-manifolds},
        date={1986},
        ISSN={0022-040X},
     journal={J. Differential Geom.},
      volume={24},
      number={3},
       pages={343\ndash 372},
         url={http://projecteuclid.org/euclid.jdg/1214440552},
      review={\MR{868975}},
}

\bib{RG18}{article}{
      author={Ramos, Vinicius Gripp~Barros},
       title={Absolute gradings on {ECH} and {H}eegaard {F}loer homology},
        date={2018},
        ISSN={1663-487X},
     journal={Quantum Topol.},
      volume={9},
      number={2},
       pages={207\ndash 228},
         url={https://doi.org/10.4171/QT/107},
      review={\MR{3812797}},
}

\bib{Sa18}{article}{
      author={Sato, Kouki},
       title={Topologically slice knots that are not smoothly slice in any
  definite 4-manifold},
        date={2018},
        ISSN={1472-2747},
     journal={Algebr. Geom. Topol.},
      volume={18},
      number={2},
       pages={827\ndash 837},
         url={https://doi.org/10.2140/agt.2018.18.827},
      review={\MR{3773740}},
}

\bib{Ta10}{article}{
      author={Taubes, Clifford~Henry},
       title={Embedded contact homology and {S}eiberg-{W}itten {F}loer
  cohomology {I}},
        date={2010},
        ISSN={1465-3060},
     journal={Geom. Topol.},
      volume={14},
      number={5},
       pages={2497\ndash 2581},
         url={https://doi.org/10.2140/gt.2010.14.2497},
      review={\MR{2746723}},
}

\bib{Thu}{article}{
      author={Thurston, William},
       title={The geometry and topology of 3-manifolds},
        date={1978},
     journal={Princeton lecture notes},
         url={http://library.msri.org/books/gt3m/},
}

\bib{Wa64}{article}{
      author={Wall, C. T.~C.},
       title={Diffeomorphisms of {$4$}-manifolds},
        date={1964},
        ISSN={0024-6107},
     journal={J. London Math. Soc.},
      volume={39},
       pages={131\ndash 140},
         url={https://doi.org/10.1112/jlms/s1-39.1.131},
      review={\MR{163323}},
}

\bib{wall}{article}{
      author={Wall, C. T.~C.},
       title={On simply-connected {$4$}-manifolds},
        date={1964},
        ISSN={0024-6107},
     journal={J. London Math. Soc.},
      volume={39},
       pages={141\ndash 149},
         url={https://doi.org/10.1112/jlms/s1-39.1.141},
      review={\MR{163324}},
}

\bib{Wa18}{article}{
      author={Watanabe, Tadayuki},
       title={Some exotic nontrivial elements of the rational homotopy groups
  of {$\mathrm{Diff}(S^4)$}},
        date={2018},
      eprint={arXiv:1812.02448},
}

\bib{Wa20}{article}{
      author={Watanabe, Tadayuki},
       title={Theta-graph and diffeomorphisms of some 4-manifolds},
        date={2020},
      eprint={arXiv:2005.09545},
}

\end{biblist}
\end{bibdiv}

\end{document}